\documentclass[11pt,reqno]{amsart}

\usepackage{CBstyle}
\usepackage[T1]{fontenc}
\usepackage[english]{babel}

\usepackage[margin=1in,footskip=0.2in]{geometry}

\title{Equidistribution of continuous low-lying pairs of horocycles via Ratner's theorem}
\author{Claire Burrin}
\address{Institute of Mathematics, University of Zurich, Switzerland}
\email{claire.burrin@math.uzh.ch}

\thanks{Thanks are due to Alex Gorodnik for organizing the workshop {\em Distribution of orbits:~Arithmetics and Dynamics} in the Swiss Alps, where Valentin Blomer's nice  talk motivated this note. Following his presentation, Amir Mohammadi and Andreas Str\"ombergsson both independently also pointed out the Ratner-based argument presented in this note. The computation that follows benefitted from separate discussions with Andreas Str\"ombergsson, and I want to thank Jens Marklof for references and comments, and the referees for their careful reading and helpful suggestions to improve the exposition of this note. The author acknowledges the support of the Swiss National Science Foundation, Grant No.~201557.}

\begin{document}

\maketitle

\begin{abstract}
We record an alternative proof of a recent joint equidistribution result of Blomer and Michel, based on Ratner's topological rigidity theorem. This approach has the advantage of extending to non-uniform lattices.
\end{abstract}

\vspace{.5cm}

Expanding closed horocycles are known to become equidistributed on any non-compact finite-area hyperbolic orbifold, such as, e.g., the low-lying horocycle $\{ \SL_2(\Z)(x+\tfrac{i}{T})\, \mid\, x\in[0,1]\}$  on the modular curve $X=\SL_2(\Z)\bk \h$ as $T\to \infty$. Equidistribution also holds for various arithmetically interesting growing discrete subsets of points along expanding closed horocycles (see, e.g., \cite{Hejhal1996,ClozelOhUllmo2001,MarklofStrombergsson2003,SarnakUbis2015,BurrinShapiraYu2022}). Most recently Blomer and Michel \cite{BlomerMichel2023} investigated the joint distribution properties of sets of the shape
\begin{align}\label{mixed}
\left\{\left(\SL_2(\Z)\left(\frac{a+i}{q}\right),\SL_2(\Z)\left(\frac{ab+i}{q}\right)\right) \mid \, a\, ({\rm mod }\, q)\right\}\subset X\times X
\end{align}
as $q\to\infty$ while $b\in (\Z/q\Z)^*$ varies with $q$ in a prescribed manner. In parallel, they consider a "continuous version" of this setup, which establishes the joint equidistribution of two low-lying horocycles of different speeds under a weak diophantine condition; see Theorem \ref{thm:BM} below. By Dirichlet's  theorem, given any real numbers $y$ and $Q>0$, there is a rational number $\tfrac{a}{q}$ so that $|y-\tfrac{a}{q}|< \tfrac{1}{qQ}$ and $q\leq Q$. In the following statement, $y$ is irrational and $Q$ is picked to depend on a growing parameter $T$ so that as $T\to\infty$ we can find a sequence $q\to\infty$ for which $y$ has a rational approximation of that order of magnitude.

\begin{Thm}[Blomer--Michel]\label{thm:BM}
Let $X=\SL_2(\Z)\bk \h$, $T>1$, $y\in[1,2]$ and write $y=a/q+O(1/(qQ))$ for positive coprime integers $a, q$ with $q\leq Q:= T^{0.99}$. Let $I\subseteq\R$ be a fixed non-empty interval. Then 
$$
\left\{\left(\SL_2(\Z)\left(x+\frac{i}{T}\right),\SL_2(\Z)\left(xy+\frac{i}{T}\right)\right) \vert x\in I\right\} \subseteq X\times X
$$
equidistributes as $T\to \infty$ for pairs $(y,T)$ with $q\to\infty$.\footnote{ In the statement of \cite[Thm 1.3]{BlomerMichel2023}, $x$ is replaced by $x/T$. The argument of proof does however refer to the low-lying horocycles stated above; see \cite[Thm 4.2]{BlomerMichel2023} and \cite[p.55]{BlomerMichel2023}.}
\end{Thm}

The argument of proof is based on the estimation of Weyl sums via a shifted convolution problem, Sato-Tate and a sieving argument, and crucially the measure classification theorem of Einsiedler--Lindenstrauss \cite[Thm 1.4]{EL19}. The purpose of this note is to record the following version of their result, using exclusively homogeneous dynamics. It would be very interesting to have such an argument for the discrete equidistribution analogue of Theorem \ref{thm:BM} concerning sets of the shape (\ref{mixed}), which is the main result of \cite{BlomerMichel2023} (see  \cite[Thm 1.1]{BlomerMichel2023}).

Write $a_t=\bsm \sqrt{t} &0\\ 0&\sqrt{t}^{-1}\esm$ for $t\in \R_{>0}$ and $u_x=\bsm 1 & x\\ 0& 1\esm$ for $x\in\R$.

\begin{Thm}\label{thm:equid}
Let $G=\SL_2(\R)$, let $\G$ be a non-uniform lattice in $G$, set $L=G\times G$, $\Lambda = \G\times \G$. 
Let $T>1$, and let $y>0$ be irrational if $\G$ is arithmetic or $y\neq1$ if $\G$ is non-arithmetic, and let $I\subseteq \R$ be a fixed non-empty interval. Then
\begin{align*}
   \left\{ \left( \G u_x a_T^{-1}, \G u_{xy} a_T^{-1}\right): x\in I\right\} \subseteq \Lambda \setminus  L
\end{align*}
equidistributes as $T\to\infty$.
\end{Thm}

\begin{proof}
Let $\pi$ be the composition of the smooth immersion $\iota_y: G\to L$, $\iota_y(g) =  (g, a_y g a_y^{-1})$ and the canonical projection $L \to \Lambda \setminus L$ so that  for any test function $\varphi \in C_b(\Lambda \bk  L)$, we have
\begin{align*}
  \int_0^1 \varphi(\pi(u_x a_{1/T}))\, dx =  \int_0^1 \varphi(\G u_x a_{1/T}, \G u_{xy}a_{1/T})\, dx,
\end{align*}
where we use the standard fact that the expansion (or contraction) of closed horocycles by the geodesic flow is algebraically expressed by $a_y u_x a_y^{-1} = u_{xy}$. A result of Shah \cite[Theorem 1.4]{Shah1996} asserts that if $\overline{\pi(G)}
=\Lambda \backslash  L$  then
\begin{align*}
 \lim_{T\to\infty} \int_0^1 \varphi(\G u_x a_{1/T}, \G u_{xy}a_{1/T})\, dx = \int_{\Lambda \setminus  L} \varphi.
\end{align*}
We may introduce a conjugation so that $L=G\times G$, $\Lambda = \G \times a_y^{-1} \G a_y$ and $\pi$ corresponds to the diagonal embedding. Let $\Delta(G)$ denote the diagonal embedding of $G$ in $L$. By Ratner's topological rigidity theorem \cite{Ratner1991}, we know that $\overline{\pi(G)}=\pi(H)$ for some closed connected subgroup $H < L$ containing $\Delta(G)$, where $H\cap \Lambda$ is a lattice in $H$.  Since $G$ is simple, $\Delta(G)$ is a maximal connected subgroup in $L$ and so the only options are that $H=L$ (in which case $\overline{\pi(G)}=\Lambda\bk L$ and equidistribution follows) or that $H=\Delta(G)$. In the latter case, the requirement that $\Lambda\cap H$ be a lattice in $H$ is equivalent to asking that $\G_y := \G \cap a_y^{-1}\G a_y$ be a lattice in $G$. This holds if and only if $a_y\in {\rm Comm}(\Gamma)$, where ${\rm Comm}(\G)$ is the commensurator of $\G$ in $G$, i.e.,
$$
{\rm Comm}(\G) = \{ \alpha \in \GL_2^+(\R) \mid \G\cap \alpha^{-1}\G \alpha \text{ has finite index in both } \G \text{ and } \alpha^{-1}\G \alpha\}.
$$

A non-uniform lattice $\Gamma<G$ is arithmetic if and only if it is commensurable to $\SL_2(\Z)$. Since ${\rm Comm}(\SL_2(\Z)) = \R^\times \cdot \GL_2^+(\Q)$, we have $a_y \in {\rm Comm}(\G)$ if and only if $y\in \Q$. If $\G$ is non-arithmetic, we claim that $a_y\in {\rm Comm}(\G)$ if and only if $y=1$. 

By Margulis' arithmeticity criterium \cite[Thm 1.16, Chapter IX]{Margulis1991} ${\rm Comm}(\G)$ contains $\G$ as a finite-index subgroup. In particular, ${\rm Comm}(\G)$ is itself a lattice in $G$ and by discreteness, its maximal unipotent subgroup is infinite cyclic. Let $u_t$ be its generator, i.e., $u_{t'}\not\in {\rm Comm}(\G)$ for all $|t'|<t$. Let $a_y\in {\rm Comm}(\G)$; up to replacing $a_y$ by its inverse, we may assume that $y>1$. Then the fact that  $a_y^{-1} u_t a_y = u_{t/y}\in {\rm Comm}(\G)$ contradicts the minimality of $t$. 
\end{proof}

\begin{Rmk}
The statement extends to finitely many factors as follows. If $y_1,\dots,y_n \in \R$ are chosen such that $a_{y_i/y_j}\not\in \mathrm{Comm}(\G)$ for $i\neq j$, then
\begin{align*}
\{ (\G u_{xy_1}a_T^{-1},\dots, \G u_{xy_n} a_T^{-1}) \mid  x\in I\} \subset \G^n\bk G^n
\end{align*}
equidistributes as $T\to\infty$. This can be compared to \cite[Thm 10.2]{Marklof2007} and \cite[Thm 4]{Strombergsson2004}. See \cite{Marklof2010} for a higher rank version, which also relies on commensurability.
\end{Rmk}

In the rest of this note, we come back to the case of $X=\SL_2(\Z)\bk \h$. By spectral expansion and Weyl's equidistribution criterion, equidistribution in Theorem \ref{thm:BM} holds if and only if 
$$
\int_0^1 \varphi\left(x+\tfrac{i}{T}, xy+ \tfrac{i}{T}\right)\, dx = o_\varphi(1)
$$ 
holds for all $\varphi=f_1\otimes f_2$, where $f_1$ and $f_2$ are non-constant Hecke eigenforms. We will show that when $y\in\Q$ (and some additional constraints on $y$) this fails, with a limit that can be expressed as a product of Hecke eigenvalues and matrix coefficients. To proceed we first recall some notation connected to the construction of Hecke operators (see, e.g., \cite{Miyake1989}). For each element $\alpha \in {\rm Comm}(\G)$ we have 
\begin{align*}
   \G \alpha  \G = \sqcup_{m=1}^M \G \alpha  h_m  
\end{align*}
where the $h_m$'s are coset representatives $\G = \sqcup_{m=1}^M \G_\alpha h_m$, with $\G_\alpha = \G \cap\alpha^{-1}\G \alpha$. For every function $f$ on $\G\bk G$ we set $(f\vert \G \alpha \G)(g) = \sum_{m=1}^M f(\alpha h_m g)$. Further, if $\alpha \in \GL_2^+(\Q)$ has integer entries, there is a uniquely determined pair of positive integers $l,m$ such that $l\mid m$ and
$ \G \alpha  \G = \G \bsm l &0 \\0 & m \esm \Gamma.
$
The Hecke operator arising from the double coset $\G \alpha \G$ is then given by 
$$
T_n f = n^{-1/2} \sum_{\substack{lm=n\\ l\mid m}} (f\vert \G\bsm l &0\\ 0&m\esm \G).
$$
We say that $f$ is a Hecke eigenform if $f$ is an eigenfunction of the Hecke operators $T_n$, $n\geq1$. We denote the corresponding Hecke eigenvalues $\lambda_f(n)$.

By abuse of notation, we view functions on $X$ as functions on $G$ that are left $\SL_2(\Z)$-invariant and right $\mathrm{SO}(2)$-invariant. We denote by $\sigma_1$ the usual sum-of-divisors function $\sigma_1(n)=\sum_{d\mid n} d$.

\begin{Prop}
Let $\G=\SL_2(\Z)$. Let $y=\tfrac{p}{q}$ with $p,q$ distinct, coprime and squarefree, and let $\varphi=f_1\otimes f_2\in C_b(X\times X)$, where $f_1, f_2$ are non-constant Hecke eigenforms. Then 
 $$
 \lim_{T\to\infty} \int_0^1 \varphi\left(x+\tfrac{i}{T}, \tfrac{xp}{q}+\tfrac{i}{T}\right)\, dx  = \frac{\lambda_{f_2}(pq)\sqrt{pq}}{\sigma_1(pq)} \scal{f_1, a_{p/q}.f_2} = \frac{\lambda_{f_1}(pq)\sqrt{pq}}{\sigma_1(pq)} \scal{ a_{q/p}.f_1,f_2},
 $$
 where $\scal{f_1,a_{p/q}.f_2}=\int_{\G\bk G} f_1(g) f_2(ga_{p/q}^{-1})d\mu(g)$ and $\mu$ is normalized so that $\mu(\G\bk G)=1$.
 \end{Prop}

\begin{proof}
Set $\G=\SL_2(\Z)$. Since $a_{p/q}\in {\rm Comm}(\G)$, the group $\G_{p/q}= a_{p/q}^{-1}\G a_{p/q}\cap \G$ is a finite index subgroup of $\G$; in particular $\G_{p/q}<G$ is a lattice. Then $\varphi\circ \iota_{p/q}\in C_b(\G_{p/q} \backslash  G)$, and the equidistribution of pieces of expanding closed horocycles on $\G_{p/q}\backslash G$  yields \cite{Strombergsson2004}
\begin{align*}
   \lim_{T\to \infty} \int_I \varphi(u_{x}a_{T}^{-1}, u_{xp/q}a_T^{-1})\, dx = \frac{1}{\mu(\G_{p/q} \backslash  G)} \int_{\G_{p/q} \backslash  G} \varphi(g, a_{p/q}ga_{p/q}^{-1}) \, d\mu(g),
\end{align*}
where the Haar measure $\mu$ is normalized so that $\mu(\G \backslash  G)=1$.

For $p,q$ distinct, coprime and squarefree, we find $\G a_{p/q}  \G = \delta_{\sqrt{pq}}^{-1} \G\bsm 1 &0 \\ 0 &pq\esm \G$, where $\delta_{x}= \bsm x & \\ & x\esm$, and $ \mu(\G_{p/q}\backslash G)=M=\sigma_1(pq)$ \cite[Lemma 4.5.6]{Miyake1989}. Let $\mathcal{F}$ be a fundamental domain for $\Gamma\bk G$; then $\sqcup_{m=1}^M h_m \mathcal{F}$, where the $h_m$'s are coset representatives of $\G_{p/q}$ in $\G$, is a fundamental domain for $\G_{p/q}\bk G$. We can now compute 
\begin{align*}
  \frac{1}{\mu(\G_{p/q}\backslash G)} \int_{\G_{p/q}\backslash G} f_1(g) f_2(a_{p/q} g a_{p/q}^{-1})\, d\mu(g) &= \frac{1}{M}\sum_{m=1}^M \int_{h_m \mathcal{F}} f_1(g) f_2(a_{p/q} g a_{p/q}^{-1})\, d\mu(g) \\
  &=
 \frac{1}{M}  \int_{\mathcal{F}} f_1(g)\sum_{m=1}^M f_2(a_{p/q} h_m g a_{p/q}^{-1})\, d\mu(g)\\
  &= \frac{\sqrt{pq}}{\sigma_1(pq)} \int_{\G \bk G} f_1(g) T_{pq} f_2(ga_{p/q}^{-1}) d\mu(g).
\end{align*}
\end{proof}

\end{document}